\documentclass[12pt]{article}
\usepackage{theorem,amsfonts}

\newtheorem{theorem}{Theorem}[section]
\newtheorem{proposition}{Proposition}[section]
\newtheorem{lemma}{Lemma}[section]
\newtheorem{corollary}{Corollary}[section]
\theorembodyfont{\upshape}
\newtheorem{definition}{Definition}[section]
\newtheorem{remark}{Remark}[section]
\newtheorem{proof}{Proof}

\newcommand{\bt}{\begin{theorem}}
\newcommand{\et}{\end{theorem}}
\newcommand{\bl}{\begin{lemma}}
\newcommand{\el}{\end{lemma}}
\newcommand{\bp}{\begin{proposition}}
\newcommand{\ep}{\end{proposition}}
\newcommand{\bo}{\begin{proof}}
\newcommand{\eo}{\end{proof}}
\newcommand{\bd}{\begin{definition}}
\newcommand{\ed}{\end{definition}}
\newcommand{\br}{\begin{remark}}
\newcommand{\er}{\end{remark}}
\newcommand{\bc}{\begin{corollary}}
\newcommand{\ec}{\end{corollary}}
\newcommand{\be}{\begin{enumerate}}
\newcommand{\ee}{\end{enumerate}}

\title{Identity excluding groups}
\author{C. R. E. Raja}
\date{}

\begin{document}
\maketitle

\let\epsi=\epsilon
\let\vepsi=\varepsilon
\let\lam=\lambda
\let\Lam=\Lambda 
\let\ap=\alpha
\let\vp=\varphi
\let\ra=\rightarrow
\let\Ra=\Rightarrow 
\let\LRa=\Leftrightarrow
\let\Llra=\Longleftrightarrow
\let\Lla=\Longleftarrow
\let\lra=\longrightarrow
\let\Lra=\Longrightarrow
\let\ba=\beta
\let\ga=\gamma
\let\Ga=\Gamma
\let\un=\upsilon

\begin{abstract}
We consider identity excluding groups.  We first show that motion
groups of totally disconnected nilpotent groups are identity
excluding.  We prove that certain class of p-adic algebraic groups
which includes algebraic groups whose solvable radical is type $R$ have 
identity excluding property.  We also prove
the convergence of averages of representations for some 
solvable groups which are not necessarily identity excluding.  
\end{abstract}

\begin{section}{Introduction}

Let $G$ be a locally compact, $\sigma$-compact 
metrizable group with a right
invariant Haar measure $m$.   Let ${\cal P}(G)$ be the space of
regular Borel probability measures on $G$.  The convolution of any two
measures $\mu$ and $\lam$ in ${\cal P}(G)$ is defined by $\mu*\lam (f) =
\int \int f(xy) d\mu (x) d\lam (y)$ for all continuous bounded functions
on $G$.  Let $L^p (G)$, $1\leq p <\infty$ be the space of all measurable
functions $f$ with $\int |f|^p < \infty$ and $L^\infty (G)$ be the space
of all (essentially) bounded measurable functions on $G$.  

For $\mu , \lam \in {\cal P}(G)$ and $x \in G$, $\mu \lam$, $x \mu$ and 
$\mu x$ denote $\mu * \lam$, $\delta _x*\mu $ and $\mu *\delta _x$ and
$\mu ^n$ denotes the $n$-th convolution power of $\mu$.  

For $\mu \in {\cal P}(G)$ and $f \in L^1 (G)$, we define the convolution
operator $\mu *f$ by $\mu *f (x) = \int f(xy) d\mu (y)$ for all $x \in G$.  

We now recall that a measure $\mu$ in ${\cal P}(G)$ is called {\it
ergodic} if ${1\over n}\sum \mu ^k *f \ra 0$ for all $f \in L^1 (G)$ with
$\int f =0$ and $\mu$ is called {\it weak mixing} if 
${1\over n} \sum |<\mu ^k*f , g>| \ra 0$ for all $g \in L^\infty (G)$ and
all $f \in L^1(G)$ with $\int f =0$.  

It is known that any weak mixing measure is ergodic and also strictly
aperiodic.  The converse of this is known as {\it weak mixing problem}.  
Lin and Wittmann show that strong convergence of $\mu ^n$-averages of 
unitary representation gives affirmative answer to weak mixing problem
(Theorem 3.1 of [LW]).  
Thus, the relationship between ergodic and weak mixing
is closely related to the representation theoretic property, called
identity excluding.  Here we investigate the class of identity excluding
groups.  We first observe that for groups whose connected
component of identity is compact, trivial representation is not weakly 
contained in any non-trivial irreducible representation implies
identity excluding.  Using this observation we prove
that (i) motion groups of totally disconnected nilpotent groups are
identity excluding and (ii) a class of p-adic algebraic groups which
includes groups whose solvable radical is type $R$ and the general affine 
group has identity excluding property.  We also show the convergence of
averages of unitary representation for splitable solvable 
Zariski-connected $p$-adic algebraic groups.
\end{section}

\begin{section}{Preliminaries}

We introduce notions and prove lemmas that are needed to prove the main
results.

\bd
A unitary representation of a locally compact group $G$ in a
Hilbert space $\cal H$ is a homomorphism $T\colon G \ra {\cal B}_u({\cal
H})$, the space of unitary operators on $\cal H$ such that for each $v
\in V$, the map $g \mapsto T(g) v$ is continuous.  
For any subgroup $H$ of
$G$ and any representation $T$ of $H$, $T|_H$ denotes the representation
of $H$ obtained by restricting $T$ to $H$.  Let $I$ be denote the trivial
irreducible representation.
\ed

Throughout this article by a group we mean a locally compact,
$\sigma$-compact metric group and by a representation we mean an unitary
representation.

\bd
We say that a representation $T$ of a group $G$ {\it weakly contains the
trivial representation} or $I\prec T$ if there exists a sequence $(v_n)$
of unit vectors in $\cal H$ such that $||T(g) v_n -v_n || \ra 0$ for all
$g\in G$: such sequence $(v_n)$ of unit vectors is called {\it
approximate fixed point}.  
\ed

\br
The standard definition of weak containment of trivial representation
requires the convergence $||T(g) v_n -v_n|| \ra 0$ to be uniform on
compact subsets of $G$.  It may be seen as follows that our definition is
equivalent to the standard definition: if there exists a sequence $(v_n)$
such that $||T(g) v_n -v_n || \ra 0$ for all $g \in G$, then for any $f
\in L^1(G)$ with $f \geq 0$ and $\int f =1$, we have $||T(f)||=1$ and
hence by 1.3, Chapter 3 of [M], $T$ weakly contains the trivial
representation according to the standard definition.
\er 

\bd
We say that a group $G$ has {\it identity excluding property} if there
is no non-trivial
irreducible representation $T$ for which there is a dense subgroup $D$ of
$G$ such that $I_D\prec T|_D$.
\ed  

\bd
For any
$\mu \in {\cal P}(G)$ and any representation $T$ of $G$, the $\mu$-average
$T_\mu$ is defined as $T_\mu (v) = \int T(g) v d\mu (g) $ for any vector
$v$.  It is easy to see that $||T_\mu ||\leq 1$.  
\ed

We now recall two non-degeneracy conditions for measures on groups which 
are necessary for weak mixing.

\bd
A $\mu \in {\cal P}(G)$ is called {\it adapted} if the closed subgroup
generated by the support of $\mu$ is $G$ and $\mu$ is called {\it
strictly aperiodic} if there is no proper closed normal subgroup a coset
of which contains the support of $\mu$.
\ed

Identity excluding groups was introduced in [JRW] and also considered in
[LW] and [Ra2].  It is shown in [LW] that nilpotent groups are identity
excluding.  Also, Lin and Witmann proved that for adapted, strictly
aperiodic measure $\mu$ on a group $G$ with identity excluding property, 

\be
\item $||T_\mu ^n|| \ra 0$ for any non-trivial irreducible unitary
representation $T$, 

\item $(T_\mu ^n)$ converges strongly for any unitary representation $T$
and 

\item in addition if $\mu$ is ergodic, then $\mu$ is weak mixing. 
\ee

We now prove some lemmas which are needed in the sequel. 
 
\bl\label{rl1}
Let $G$ be a group and $H$ be a closed abelian normal subgroup of $G$ such
that $G/C(H)$ is finite where $C(H)$ is the centralizer of $H$.  Suppose
$T$ is an irreducible representation of $G$ such that $I\prec T$.  Then 
$T$ is trivial on $H$.

\el

\bo
Let $T$ be an irreducible representation of $G$ in a Hilbert space 
${\cal H}$.  Suppose $I\prec T$.  
For any character $\chi$  of $H$, define 
$V_\chi = \{ v\in {\cal H}\mid T(h) v = \chi(h) v {~~\rm for ~all~~}
h \in H \}$.  Then it is easy to see that $V_\chi$ is a $C(H)$-invariant
closed subspace of $\cal H$.  Let $\chi $ be such that $V_\chi$ is
a non-trivial subspace.  Let $g_0C(H), g_1C(H), \cdots g_k C(H)$ be a
system of coset representative of $C(H)$ in $G$ with $g _0 = e$.  Now for
each $0\leq i \leq k$, define $\chi _i (h) = \chi (g_ih g_i^{-1})$ for all
$ h \in H$.  Then each $\chi _i$ is a character of $H$ and $T(g_i)V_\chi =
V_{\chi _i}=V_i$, say.   Then the orthogonal sum $\oplus V_i$ is a
$G$-invariant closed subspace of $\cal H$.  Thus, ${\cal H}=\oplus V_i$.  
Since $I\prec T$, there exists a $i\geq 0$ such that $I\prec T_i$ where
$T_i$ is the representation of $H$ in $V_i$ given by $\chi _i$.  This
shows that $\chi _i$ is trivial and hence $\chi$ is trivial.  Thus, $T(H)$ 
is trivial.
\eo

\bl\label{rl2}
Let $G$ be a group whose connected component of identity is
compact and $T$ be an unitary representation of $G$ in a Hilbert space 
$\cal H$.  Suppose there is a
dense subgroup $H$ of $G$ and a sequence $(v_n)$ of unit vectors in 
$\cal H$ such that $||T(g) v_n -v_n || \ra 0$ for all $g \in H$.  Then 
$||T(g)v_n -v_n || \ra 0$ for all $g \in G$.  In particular, $I\not \prec
T$ for any non-trivial irreducible representation of $G$ implies $G$ is
identity excluding.
\el

\bo
Let $G^0$ be the connected component of identity in $G$.  Then $G/G^0$ is
totally disconnected and has compact open subgroups, see Theorem II.7.7 of
[HR].  Since $G^0$ is compact, $G$ itself has compact open subgroups.  
Let $K$ be a compact open subgroup of $G$.  Let 
${\cal H}_K=\{v\in {\cal H}\mid T(K)v=v \}$.  Now for
each $n \geq 1$, there exists a $a_n \in {\cal H}_K$ and a $b_n \in {\cal
H}_K^\perp$ such that $v_n = a_n +b_n$.  It can be easily seen that
${\cal H}_K^\perp$ is also $T(K)$-invariant.  This implies that for 
$g\in H\cap K$, $||T(g)v_n -v_n|| = ||T(g)b_n-b_n||\ra 0$.  
Since $H$ is dense in $G$, $H\cap K$ is dense in $K$.  We may assume that
$H$ is a Borel subgroup of $G$.  Thus, $H\cap K$ is a Borel subgroup of
$K$.  Since $K$ is metrizable, 
we can choose a countable dense subset $E$ of $K$ such that $E\subset 
H\cap K$.  Let $\mu = \sum _{x\in E} r_x \delta _x$ where $r_x>0$
for all $x\in E$ and $\sum r_x =1$.  Then $\mu$ is an adapted, 
strictly aperiodic probability measure on $K$ such that 
$\mu (H\cap K)=1$.  Let $\rho$ be the representation of $K$ in ${\cal
H}_K^\perp$ such that $\rho (g)$ is the restriction of $T(g)$ to
${\cal H}_K^\perp$ for any $g \in K$.  Then 
$||(\rho _\mu  ^k -I)b_n||\ra 0$ as $n \ra \infty$.  
Suppose $\rho _\mu ^k -I$ is not invertible for any $k \geq 1$.  
Then $||\rho _\mu ^k || =1$ for all $k \geq 1$.  
Since $K$ is a SIN-group by Theorem 2.11 of [LW], $I\prec \rho$.  Since
$K$ has property $(T)$, ${\cal H}_K^\perp$ has a non-trivial $T(K)$-fixed
vector.  This is a contradiction.  This shows that $(\rho _\mu ^k -I)$ is
invertible for some $k \geq 1$.  Thus, $b_n \ra 0$.  Now for any 
$g \in K$, $||T(g)v_n -v_n ||= ||T(g) b_n -b_n|| \leq 2||b_n ||
\ra 0$.  Thus, $\{g \in G \mid ||T(g) v_n -v_n|| \ra 0 \}$ is a
dense subgroup containing an open subgroup $K$ and hence 
$||T(g)v_n -v_n ||\ra 0 $ for all $g \in G$.
\eo
\end{section}

\begin{section}{Motion groups}

We first prove the identity excluding for motion groups of totally
disconnected nilpotent groups. 

\bt\label{t1}
Let $G$ be a split compact extension of a totally disconnected nilpotent
group.  Then $G$ has identity excluding property.
\et

\bo
Let $N$ be a totally disconnected nilpotent normal subgroup of $G$ such
that $G/N$ is compact.  Since $G$ is a split compact extension of $N$, 
it is easy to see that the
connected component of $G$ is compact.  In view of Lemma \ref{rl2},
it is enough to show that $I\not \prec T$ for any non-trivial irreducible
representation.  

Let $T$ be an irreducible representation of $G$ such that $I\prec T$.  
Let $Z$ be the center of $N$.  We now prove that $T (Z)$ is trivial.  Let
$K$ be a compact open subgroup of $Z$ and define $L = \cap _{g \in
G}gKg^{-1}$.  Then $L$ is a compact normal of subgroup of $G$.  Since
$I\prec T$, $T(L)$ has a non-trivial fixed point.  Since $L$ is a normal 
subgroup of $G$, space of fixed points for $T(L)$ is a non-trivial
$T(G)$-invariant subspace.  Thus, $T(L)$ is trivial.  
Since $K$ is a open subgroup of $Z$, $N(K)$, the normalizer of $K$ is an
open subgroup of $G$ containing $N$.  Since $G/N$ is compact, $N(K)$ is a
open subgroup of finite index in $G$.  Thus, $L$ is a finite intersection
of conjugates of $K$.  Since $Z$ is a normal subgroup of $G$, $L$ is
a open subgroup of $Z$.  Now by replacing $G$ by $G/L$, we may assume that
$Z$ is discrete.  Let $x \in Z$.  There exists a finitely generated normal
subgroup $F$ of $G$ such
that $x \in F \subset Z$ (we may choose $F$ to be the group generated by
the finite set $\{g xg^{-1} \mid g \in G \}$).  It is easy to see that
Aut~$(F)$ is a discrete group.  Since $H\subset Z$, the center of $N$ and
$G/N$ is compact, the centralizer of $F$ is finite.  By Lemma \ref{rl1},
$T(F)$ is trivial.  In particular, $T(x)$ is trivial.  Since $x \in Z$ is
arbitrary, $T(Z)$ is trivial.  

Since $N$ is nilpotent, by using the
central series of $N$, we conclude that $T(N)$ is trivial.  Since any
compact group is identity excluding, $T(G)$ itself is trivial.  Thus,
$G$ is identity excluding.
\eo

We next prove that finite extension of nilpotent groups are 
identity excluding. 

\bt\label{fe}
Let $G$ be a finite extension of a nilpotent group.  Then $G$ is 
identity excluding.
\et

\bo
Let $N$ be a nilpotent subgroup of $G$ such that $G/N$ is finite.  Then
replacing $N$ by a subgroup of $N$, we may assume that $N$ is a normal
subgroup of $G$ (the subgroup may be chosen to be $\cap xNx^{-1}$).  

Let $T$ be an irreducible representation of $G$ in a Hilbert space $\cal
H$.  Suppose there exists a dense set $D$ of $G$ such that $I\prec
T|_D$.  Since $G/N$ is finite, $D\cap N$ is dense in $N$.  Let $N_0 = N$
and $N_i =[N, N_{i-1}]$ for all $i >0$.  Then $N_k$ is contained in the
center of $N$ for some $k \geq 0$.  Also, $D\cap N_k$ is dense in $N_k$.  

We now claim that $T$ is trivial.  For any character $\chi$ of $N_k$,
we define ${\cal H}_\chi = \{ v\in {\cal H} \mid T(g) v = \chi (g) v
~~{\rm for ~~  all~~} g \in N_k \}$.   Let $C$ be the centralizer of $N_k$
in $G$.  Then $C$ is a normal subgroup of finite index in $G$.  Also,
${\cal H}_\chi$ is $C$-invariant closed subspace of $\cal H$.  Let $\chi $
be such that ${\cal H}_\chi$ is a non-trivial subspace.  Let $g_0C, g_1,
\cdots g_k C$ be a system of coset representative of $C$ in $G$ with $g _0
= e$.  Now for each $0\leq i \leq k$, define $\chi _i (x) = \chi (g_ix
g_i^{-1})$ for all $x\in N_k$.  Then each $\chi _i$ is a character of
$N_k$ and $g_i{\cal H}_\chi = {\cal H}_{\chi _i}={\cal H}_i$, say.  
Then the orthogonal sum $\oplus {\cal H}_i$ is a $T(G)$-invariant closed
subspace of $\cal H$.  Thus, ${\cal H}=\oplus {\cal H}_i$.  
Since $I\prec T|_D$,
there exists a $i\geq 0$ such that $I\prec T_i|_{D\cap N_k}$ where
$T_i$ is the representation of $N_k$ in ${\cal H}_i$ 
given by $\chi _i$.  This shows that $\chi _i$ is trivial and hence $\chi$
is trivial.  Thus, $T(N_k)$ is trivial.  Using the central series of $N$,
we may prove that $T(N)$ is trivial.  Since finite groups are identity
excluding, $T$ itself is trivial.  
\eo

We now deduce the following for groups of polynomial growth.  We recall
that a compactly generated group $G$ (with Haar measure $m$) is said to be
of {\it polynomial growth} if there exists an integer $l >0$ and a
constant $c>0$ such that $m (U^n) \leq cn^l$ for all $n \geq 1$ where 
$U$ is a compact neighborhood of identity generating $G$: see [Gr], [L1]
for results on polynomial growth.    

\bc\label{cl}
Let $G$ be a compactly generated group of polynomial growth.  Suppose the
connected component of identity is compact.  Then $G$ has identity 
excluding property. 
\ec

\bo
Let $G$ be a compactly generated group of polynomial
growth.  By Theorem 2 of [L1], there exists a compact normal subgroup $K$
of $G$ such that $G/K$ is a Lie group of polynomial
growth.  Suppose the connected component of identity in $G$ is compact.  
Then $G$ has a compact open subgroup and hence $G/K$ also 
has a compact open subgroup.  Since $G/K$ is a Lie group, this implies
that the connected component of identity is a compact open normal 
subgroup.  Thus, $G$ has a compact open normal subgroup, 
let it be $L$.  So for any representation $T$, the space of $T(L)$-fixed
points are invariant under $G$ and hence 
it is enough to show that $G/L$ has identity
excluding property.  Since $G/L$ is a discrete group of polynomial growth, 
by Gromov's theorem in [Gr], we get that $G/L$ is a finite extension of a 
nilpotent group.  Now the result follows from Theorem \ref{fe}.
\eo

\end{section}

\begin{section}{$p$-adic algebraic groups}

We first make the following observation: 

\bp\label{p2}
Let $T$ be an unitary representation of a group $G$.  Suppose for some 
$f\in L^1(G)$ with $f\geq 0$, $T(f)$ is a compact operator.  Then 
$I\prec T$ implies $T$ has a non-trivial fixed point. 
\ep

\bo
Let $f \in L^1 (G)$ such that $f \geq 0$ and $T(f)$ is a compact
operator.  Then we may assume that $\int f =1$.  
Suppose $I\prec T$.  Then there exists a sequence $(v_n)$ of unit
vectors such that $||T(g)v_n -v_n||\ra 0$ for all $g \in G$.  This
implies that $||v_n-T(f)v_n|| \ra 0$.  Since $T(f)$ is a compact operator, 
$(T(f)v_n )$ has a convergent subsequence.  By passing to a subsequence we
may assume that $T(f) v_n \ra v$.  Since $||T(f)v_n -v_n|| \ra 0$, 
$v_n \ra v$.  This implies that $T(G)v=v$ and $||v||=1$.  Thus, $T$ has a 
non-trivial fixed point.  
\eo

Thus, the above result shows that totally disconnected CCR-groups have 
identity excluding property (see [D] for details on CCR-groups).  It is
known that only algebraic groups
that are CCR are the direct products of a semisimple group and a group of
type $R$ ([Li] and [Pu]).  Here we prove that $p$-adic
algebraic groups which are semi-direct product of semi-simple
groups and groups of type $R$ have identity excluding property, that
is any $p$-adic algebraic group has identity excluding property if 
the solvable radical is type $R$.  By an algebraic group over a local
field $\mathbb K$ of characteristic zero, we mean the group of $\mathbb
K$-points of an algebraic group defined over $\mathbb K$.  

\bd
We say that a finite-dimensional vector space $V$ over a local field of 
characteristic zero is of {\it type $R_\Ga$} where $\Ga$ is a group of
automorphisms of $V$, if the eigenvalues of each element of $\Ga$ are of
absolute value one.  We say that a Lie group $G$ is of {\it type $R_\Ga $}
where $\Ga$ is a group of Lie automorphisms of $G$, if the Lie algebra of 
$G$ is of type $R_\Ga$ and a Lie group $G$ is said to be of {\it type
$R$} if $G$ is of type $R_{{\rm Ad}~(G)}$ where Ad is the 
adjoint representation of $G$: see [J] for results on type $R$ real Lie
groups and [Ra1] for results on $p$-adic Lie groups of type $R$.
\ed

\bt\label{mt}
Let $G$ be any p-adic algebraic group.  Let $U$ be the unipotent radical
of $G$.  Let $U_0=U$, $U_i = [U, U_{i-1}]$ for $i>0$.  For any
$i >0$ and for any $G$-invariant subspace $W$ of $U_i/U_{i+1}$ 
define $\phi _{i,W} \colon G \ra GL(W)$ by $\phi _{i, W}(g) xU_{i+1} =
gxg^{-1}U_{i+1}$ for all $xU_{i+1} \in W\subset U_i/U_{i+1}$.  
Suppose for each $i>0$ and $W$ as above, either $\phi _{i, W}(G)$ is 
non-amenable or $W$ is of type $R_{\phi _{i, W}(G)}$.  Then $G$ has
identity excluding property.
\et

\br\label{1}
Suppose $G$ is an algebraic subgroup of $GL(W)$, for some
finite-dimensional vector space $W$.  Suppose $G$ is reductive.  Then $G$
is amenable implies that $G$ has no non-trivial non-compact simple
factors.  Thus, for amenable $G$, $W$ is of type $R_G$ if and only if the
center of $G$ has no split torus, that is anisotropic.
\er

\bo
Let $S$ be a reductive Levy subgroup of $G$ and $U$ be the unipotent
radical of $G$.  Then $G$ is the semidirect product of $S$ and $U$.  Let
$Z$ be the center of $U$.  Then $Z$ contains $U_k$ where $k \geq0$ is such
that $U_k \not = (e)$ but $U_{k+1} =(e)$.  Let $V= U_k$.  

We prove the result by induction on dimension of $U$.  It may be noted
that $G/W$ satisfies the hypothesis for any
irreducible $G$-invariant subspace $W$ of $V$.  
Suppose ${\rm dim}~ (U)=0$.  Then $G$ is a reductive group.  
Since $G$ is a CCR-group, $G$ has identity excluding property.  

If ${\rm dim}(U) >0$. Then $V$ is a non-trivial normal subgroup of $G$.  
Suppose $V$ contains a irreducible $G$-subspace $W$ such that the image 
of $G$ in $GL(W)$ is non-amenable.  
It is easy to see that the action of $G$ on $\hat W$ is also
irreducible.  By 4.15 and 5.15 of [S], we get that $(G, W)$ has strong
relative property
$(T)$ (see [S] for details).  Let $T$ be a irreducible representation of
$G$ such that $I\prec T$.  Define a representation $T_1$ of the semidirect
product of $G$ and $W$ by $T_1((g, w)) =
T(gw)$.  Then $I\prec T_1|_G$ and hence $T_1(W)$ has a non-trivial fixed
point.  So, $T(W)$ has a non-trivial fixed point.  Since $W$ is a
normal subgroup, $T(W)$ is trivial.  Now by induction
hypothesis $T(G/W)$ is trivial and hence $T(G)$ is trivial.

Suppose, for any irreducible $G$-subspace $W$ of $V$, the image of $G$ in
$GL(W)$ is amenable.  Then by assumption $W$ is of type 
$R_{\phi _{i, W}(G)}$.  This implies that any split torus of $G$, 
centralizes $W$ and since $U$ centralizes $W$, $\phi _{i, W}(G)$ is 
compact.  Thus, $\{ gxg^{-1} \mid g\in G \} $ is compact in $V$ for any $x
\in V$.  

Suppose $G$ is a $p$-adic algebraic group, then $V$ is a increasing union
of compact normal subgroups of $G$.  Let $(M_i)$ be a sequence of
compact normal subgroups of $G$ such that $V= \cup M_i$ and $M_i \subset 
M_{i+1}$ for all $i >0$.  Let $T$ be an irreducible representation of $G$
such that $I\prec T$.  Then for each $i >0$, $I\prec T|_{M_i}$.  For
each $i >0$, since $M_i$ has property $(T)$, $T(M_i)$ has a non-trivial
fixed point.  For each $i >0$, since $M_i$ is normal in $G$, the space
of fixed points of $T(M_i)$ is a non-trivial $T(G)$-invariant closed
subspace.  Thus, $T(M_i)$ is trivial for all $i >0$.  This shows
that $T(V)$ is trivial.  By applying induction hypothesis to $G/V$, we
get that $T$ is trivial.  Thus, by Lemma \ref{rl2}, $G$ has identity
excluding property.
\eo

\bc
Suppose $G$ is a $p$-adic algebraic group whose solvable radical is of
type $R$.  Then $G$ satisfies the hypothesis of Theorem \ref{mt} and $G$
is identity excluding.  
\ec

\bo
Suppose for some $i >0$, there is an irreducible normal subgroup $W$ of 
$U_i/U_{i+1}$ such that the image of $G$ in $GL(W)$ is amenable.  Then any
semisimple Levy subgroup of $G$ has only compact orbits in $W$.  Since
the solvable radical is of type $R$, $Gx$ is compact for any $x \in
W$.  Thus, $G$ satisfies the hypothesis of Theorem \ref{mt}.  
\eo

The following is easy to verify.
\bc
The semidirect product of $GL_n({\mathbb Q}_p)$ and
${\mathbb Q}_p^n$ (for $n >1$) verifies the hypothesis of Theorem \ref{mt}
but its solvable radical is not of type $R$. 
\ec

We now prove a partial converse to Theorem \ref{mt}.

\bt
Let $G$ be a $p$-adic algebraic group and $U$ be the unipotent radical of
$G$.  Suppose the solvable radical of $G$ is not of type $R$ and $U$ is of
type $R_L$ where $L$ is a semisimple Levy subgroup of $G$.  Then $G$ is
not identity excluding.
\et

\bo
It is enough to prove the result for a
quotient group of $G$.  Let $S$ be the solvable radical of $G$.  Since $S$
is not of type $R$, let $U'= [U,U]$, then $S/U'$ is also not of type
$R$.  Thus, replacing $G$ by $G/U'$, we may assume that unipotent radical
is abelian.  Since $S$ is not of type $R$, there exists an element $x \in
S$ such that the subspace $\{v \in U \mid x^n vx^{-n} \ra e \}$ is
non-empty and it is $G$-invariant.  Thus, again replacing $G$ by a
quotient of $G$, we may assume that the split torus $A_s$ of $S$ is
atleast one-dimensional and for any $x \in A_s \setminus K$, $x$ or 
$x^{-1}$ contracts $U$ for some compact subgroup $K$ of $A_s$.  

Now, let $\chi$ be a non-trivial character of $U$.  Let $H$ be the
stabilizer of $\chi$ in the reductive Levy part of $G$.  Define $\rho$
on $HU$ by $\rho (hu) = \chi (u)$.  Then $\rho$ defines an irreducible
unitary representation of $HU$.  Now by Mackey's normal subgroup analysis 
(see [Ma1] and [Ma2]), the induced representation $T$ of $G$ from $\rho$
is irreducible.  We now claim that $I\prec T$.  Let $G_1 = A_sHU$.  By 
assumption, the non-compact simple factors of $L$ centralizes $U$.  This
implies that $G/G_1$ is compact.  Let $T_1$ be the induced representation
of $G_1$ from $\rho$.  Then $T$ is the induced representation from $T_1$,
by Chapter III, 1.11 of [M], it is enough to show that $I\prec T_1$.  As
in example of [JRW], we can prove that $I\prec T_1$.
\eo

\end{section}

\begin{section}{Convergence of representation averages for some solvable
algebraic groups}

Lin and Witmann show that the convergence of $(T_\mu ^n)$ in the strong
operator topology 
implies that any strictly aperiodic ergodic measure is weak mixing.  
We now consider convergence of $(T_\mu ^n)$ for measures on
any solvable groups.  We have proved that solvable $p$-adic algebraic
groups are identity excluding if and only if it is of type $R$, here we 
prove that $\mu ^n$-averages of representations on split solvable
$p$-adic algebraic groups (which are not necessarily identity
excluding) are strongly
convergent.  We recall that a solvable $p$-adic algebraic group is called
{\it split} if maximal torus is splitting.  

\bt\label{pw1}
Let $G$ be a split solvable Zariski-connected $p$-adic algebraic group 
which is not of type $R$.  
Let $\mu$ be an adapted and strictly aperiodic probability measure on 
$G$.  Let $T$ be a representation of $G$.  Then $(T_\mu ^n)$ converges in
the strong topology. 
\et

\bo
In view of Theorem 2.2 of [LW], we may assume that $T$ is
irreducible.  Let $U$ be the unipotent radical of $G$.  Let $A$ be
a maximal of torus of $G$.  Since $G$ is split, $A$ is a
split torus.  Since $G$ is not of type $R$, both $U$ and $A$ are
non-trivial.  Then the center of $U$ is non-trivial.  Let $Z$ be the
center of $U$.  We may assume that $T$ has an approximate fixed
point, otherwise $T_\mu ^n $ converges.  Then $T$ is trivial on the 
center of $G$, so we may assume that $G$ has trivial center.  

Let $\hat Z$ be the dual of $Z$.  We now claim that the stabilizer of a
non-trivial point in $\hat Z$ is a proper subgroup of $G$.  Suppose there
is a non-trivial point $\chi$ in $\hat Z$ whose stabilizer is the whole
group.  Then there exists a one-dimensional subspace $W$ of $Z$, which is
a normal subgroup of $G$.  Then $W$ is either type $R_G$ or
$W$ is contracted by an element of $A$.  Suppose $W$ is type
$R_G$.  Then since $G$ is a split solvable group, $G$ action on $W$ is
trivial.  Since center of $G$ is trivial, this is a contradiction.  So, we
may assume that there is an element $g$ of $A$ contracting $W$.  Then 
$\chi (x) = g^{-n} \cdot \chi  (x) = \chi (g^n \cdot x) \ra 1$ for
any $x \in Z$.  This is a contradiction.  Thus, the stabilizer of any
non-trivial point of $\hat Z$ is a proper subgroup of $G$.  

By Mackey's theorem (see 13.3, Theorem 1 of [Kr]), 
there exists a proper subgroup $H$
containing $U$ and an irreducible representation $\rho$ of $H$ in a
Hilbert space $E$ such that the induced representation of $G$ from
$\rho$ is $T$ (up to equivalence).  

We now claim that $G$ is a semidirect product of an abelian group and $H$.  
Let $A_0 = A\cap H$.  Then $A_0$ is subtorus of $A$.  Since $A$ splits in
$G$, there exists a subtorus $A_1$ of $A$ such that $A$ is
the direct product $A_0 \times A_1$.  Since $A_0$ is a proper split torus,
$A_1$ is a non-trivial split torus. 

Hence, we identify $G/H$ with $A_1$.  Let
$m$ be Haar measure on $A_1$.  Then under the
identification, $m$ is a $G$-invariant measure on $G/H$.  

Let $L^2 (G, H, \rho )$ be the space of all measurable functions $f \colon
G \ra E$ such that 
\be
\item [(i)] $f(xh) = \rho (h) f(x)$ for all $x\in G$ and $h \in H$,

\item [(ii)] $\int _{A_1}||f(x)||^2 dm(x) < \infty$. 
\ee
Since the induced representation from $\rho $ is $T$, $T$ is defined on
$L^2(G, H,\rho )$ by $T(g) f (x) = f(gx)$ for all $f \in L^2(G,H,\rho )$
and all $x, g\in G$ (see Chapter I, 5.2 of [M]).    

Now for any $g \in G$, there exists unique $s \in A_1$ and unique $h \in
H$ such that $g = sh$.  Thus, for $f\in L^2(G,H, \rho)$, we have 

$$\begin{array}{cl}
||T_\mu  f|| ^2& \cr 
= & |\int _{A_1}<T_\mu f(x), T_\mu f (x)> dm (x)| \cr

\leq & \int _{A_1} \int _G \int _G |<f(g_1x), f(g_2x)>| d\mu (g_1) 
d\mu (g_2) dm (x) \cr

\leq & \int _{A_1} \int _G \int _G |<\rho (x^{-1}h_1x) f(s_1x), 
\rho (x^{-1}h_2 x) f (s_2x)>| d\mu d\mu dm \cr

& {\rm where}~~ g_1 = s_1h_1 ~~{\rm and }~~ g_2 = s_2 h_2 \cr

\leq & \int _{A_1} \int _G\int _G ||f(s_1x)|| ||f(s_2x)||d\mu d\mu dm. 
\end{array}$$

Let $\lam$ be the image of $\mu$ in $A_1$ under the canonical map 
$g \mapsto s$ and $R$ be the regular representation of $A_1$ in 
$L^2 (A_1)$.  Then from the above calculations we get that 
$||T_\mu ^n f|| \leq ||R_\lam ^n F||$ where 
$F(x) = ||f(x)||$ for all $x \in A_1$.  Since nilpotent
groups are identity excluding and $A_1$ is non-compact, we get that 
$||T_\mu ^n f|| \ra 0$. 
\eo

We now consider any split solvable group.  

\bt\label{pw2}
Let $G$ be any split solvable Zariski-connected $p$-adic algebraic group
and $\mu$ be an adapted and strictly aperiodic probability measures on 
$G$.  Suppose $T$ is a representation of $G$.  Then $(T_\mu ^n )$
converges strongly.
\et

\bo
Suppose $G$ is of type $R$.  Then $G$ is a nilpotent group and hence it is
identity excluding.  Suppose $G$ is not type $R$, then strong convergence
of $(T_\mu ^n )$ follows from Theorem \ref{pw1}
\eo

\end{section}

\vskip 0.25in

\noindent {C. Robinson Edward Raja, \\
Indian Statistical Institute,\\
Statistics and Mathematics Unit,\\
8th Mile Mysore Road,\\
R. V. College Post,\\
Bangalore - 560 059.\\
India.}

\noindent {creraja@isibang.ac.in\\
raja$\_$robinson@hotmail.com}

\end{document}